\documentclass[11pt]{article}
\usepackage{graphicx}
\usepackage{amsmath,latexsym,amsbsy,amssymb, color}
\usepackage{graphicx}
\def\bpm{\begin{pmatrix}}
\def\bbm{\begin{bmatrix}}
\def\ebm{\end{bmatrix}}
\def\epm{\end{pmatrix}}

\def\argmin{\hbox{arg}\!\min}
\def\ds{\displaystyle}
\def\forall{\hbox{for all}~}
\def\L{{\bf L}}

\def\bfv{{\bf v}}
\def\bfw{{\bf w}}

\def\bfb{{\bf b}}

\def\bfv{{\bf v}}
\def\ve{\varepsilon}
\def\n{\noindent}

\def\R{{\mathbb R}} 

\def\vs{\vskip 2em}

\def\v{\vskip 1em}

\def\M{{\cal M}}

\def\G{{\cal G}}
\def\C{{\cal C}}

\def\ov{\overline}

\def\Tilde{\widetilde}

\def\bega{\begin{array}}
\def\enda{\end{array}}
\def\begi{\begin{itemize}}
\def\endi{\end{itemize}}

\def\bel{\begin{equation}\label}
\def\eeq{\end{equation}}
\def\sqr#1#2{\vbox{\hrule height .#2pt
\hbox{\vrule width .#2pt height #1pt \kern #1pt
\vrule width .#2pt}\hrule height .#2pt }}
\def\square{\sqr74}
\def\endproof{\hphantom{MM}\hfill\llap{$\square$}\goodbreak}

\newtheorem{theorem}{Theorem}[section]
\newtheorem{corollary}{Corollary}[section]
\newtheorem{lemma}{Lemma}[section]

\newtheorem{remark}{Remark}[section]
\newtheorem{definition}{Definition}[section]

\begin{document}
\title{\bf  Generic Properties of Conjugate Points in  Optimal Control Problems}\vs
	\author{Alberto Bressan$^{(1)}$, Marco Mazzola$^{(2)}$, and Khai T. Nguyen$^{(3)}$\\
\\
		{\small $^{(1)}$ Department of Mathematics, Penn State University}\\
		{\small $^{(2)}$ IMJ-PRG, CNRS, Sorbonne Universit\'e,}\\
		{\small $^{(3)}$  Department of Mathematics, North Carolina State University.}
		\\
		\quad\\
		{\small e-mails:  axb62@psu.edu,~~ marco.mazzola@imj-prg.fr,~~ khai@math.ncsu.edu}}
\maketitle
\begin{abstract} 
The first part of the paper studies a class of optimal control problems in Bolza form, 
where the dynamics is linear w.r.t.~the control function.
A necessary condition is derived, for the optimality of a trajectory which starts at a conjugate point. The second part is concerned with a classical problem in the
Calculus of Variations, with free terminal point.   For a generic terminal cost $\psi\in \C^4(\R^n)$, applying the previous necessary condition
we show that the set of conjugate points is contained in the image of an 
$(n-2)$-dimensional manifold, and 
has locally bounded $(n-2)$-dimensional Hausdorff measure.
		
\quad\\
		\quad\\
		{\footnotesize
		\medskip
		
		\n {\bf AMS Mathematics Subject Classification.} 49K05, 49L12.
		}
	\end{abstract}

\maketitle

\section{Introduction}
\label{sec:1}
\setcounter{equation}{0}

Conjugate points play a key role in the study of necessary conditions, 
for problems in the Calculus of Variations and optimal control \cite{CaFr96, Cesari, Clarke, ZZ}.
The present paper intends to be a contribution to the analysis 
of conjugate points, from the point 
of view of generic theory.   Given a family of optimal control problems,
with various terminal costs, we seek properties of the set of conjugate points which are true for nearly all terminal costs $\psi\in \C^4(\R^n)$.  Here ``nearly all" is meant in the topological sense of Baire category: these properties should be true on a $\G_\delta$ set, i.e., on the
intersection of countably many open dense subsets.
As usual, $\C^k(\R^n)$ denotes the Banach space of all bounded functions with bounded, continuous  partial derivatives up to order $k$, see for example  \cite{BLN, Evans}.

Our basic setting is as follows. Consider an optimal control problem of the form
\bel{oc1}
\hbox{minimize:}\qquad J^{\tau,y}[u]~\doteq~\ds\int_\tau^T L\bigl( x(t), u(t)\bigr)\, dt + \psi\bigl(x(T)\bigr),
\eeq
where $t\mapsto x(t)\in \R^n$ is the solution to the Cauchy problem
with dynamics linear w.r.t.~the control:
\bel{oc2} \dot x(t)~=~f\bigl(x(t), u(t)\bigr)~=~f_0(x(t))+\sum_{i=1}^mf_i(x(t))\, u_i(t),\eeq
and initial data 
\bel{idy} x(\tau)= y.\eeq
Here and in the sequel, the upper dot denotes a derivative w.r.t.~time.
In (\ref{oc1}), the minimum cost is sought among all measurable functions  
$u:[\tau,T]\mapsto\R^m$.  
For  $(\tau, y)\in [0,T]\times\R^n$, the associated value function $V$ is defined as
\bel{V}
V(\tau,y)~\doteq~\inf_{u(\cdot)} J^{\tau,y}[u]\,.
\eeq
 To fix ideas, we shall consider a couple $(f,L)$ satisfying the following hypotheses.
\begi
\item[{\bf (A1)}] {\it In (\ref{oc2})
the vector fields $f_i$, $i=0,\ldots,m$, are  three times continuously differentiable and satisfy the sublinear growth condition
 \bel{sublin}
 \bigl|f_i(x)\bigr|~\leq~c_1\, \bigl(|x|+1\bigr)
 \eeq
 for some constant $c_1>0$ and all $ x\in\R^n$.
\item [{\bf (A2)}] The running cost  $L:\R^n\times\R^m\mapsto\R$ is  three times continuously differentiable and 
uniformly convex w.r.t.~$u$. Namely, for some $\delta_L>0$, 
the $m\times m$ matrix of second derivatives w.r.t.~$u$ satisfies
\bel{u-convex-L}
L_{uu}(x,u)-\delta_L\cdot \mathbb{I}_m ~\ge ~0 \qquad\qquad\forall x,u.
\eeq
Here $ \mathbb{I}_m$ denotes the $m\times m$ identity matrix.
}
\endi
The Pontryagin necessary conditions \cite{BPi, Cesari, FR} take the form
\bel{PMP1}\left\{
\bega{rl} \dot x&=~f\bigl(x,  u(x,p)\bigr),\\[3mm]
\dot p&=~- p\cdot  f_x\bigl(x,  u(x,p)\bigr)- L_x(x,u(x,p)),\enda
\right.\eeq
where $u(x,p)$ is determined as the pointwise minimizer
\bel{PMP2} u(x,p)~=~\argmin_{\omega\in \R^m} \Big\{ L(x,\omega) + p\cdot f(x,\omega)\Big\}.\eeq

The assumptions in {\bf (A1)}-{\bf (A2)} guarantee that the minimizer in (\ref{PMP2})  is unique and  solves
\bel{op1}p\cdot 
 f_u(x,\omega)+L_u(x,\omega)~=~0.
\eeq
Therefore the map $(x,p)\mapsto u(x,p)$ is well defined and continuously differentiable, and  the system of ODEs (\ref{PMP1}) has continuously differentiable right hand side. 
In particular,  for any $z\in\R^n$, the system (\ref{PMP1}) with  terminal conditions
\bel{td3} x(T)~=~z,\qquad p(T)~=~\nabla\psi(z),\eeq
admits a unique  solution  $t\mapsto (x,p) (t,z)$ defined on $[0,T]$. In turn, this uniquely determines the 
control 
\bel{uxp}t~\mapsto~ u(t,z)~=~u\bigl( x(t,z), \, p(t,z)\bigr).\eeq
In the following we mainly focus on the case $\tau=0$. 
\begin{definition}\label{def:1}
Given an initial point $\bar x\in \R^n$, we say that 
a control $u^*:[0,T]\mapsto \R^m$ is a {\bf weak local minimizer} of the cost functional
\bel{Jbx} J^{\bar x}[u]~\doteq~\int_0^T L\bigl( x(t), u(t)\bigr)\, dt + \psi\bigl(x(T)\bigr),
\eeq
subject to
\bel{dyn2}\dot x~=~f(x, u),\qquad\qquad 
x(0)~=~\bar x,\eeq
if there exists $\delta>0$ such that
$
J^{\bar x}[u^*]\leq J^{\bar x}[u]$ for every measurable control $u(\cdot)$ such that
$\|u-u^*\|_{\L^\infty} <\delta$.
\end{definition}

Consider again the maps 
$$z\,\mapsto\, x(\cdot, z),\qquad z\,\mapsto\, p(\cdot, z),\qquad z\,\mapsto\, u(\cdot, z)
$$ as  in (\ref{uxp}), obtained by solving the backward Cauchy problem
(\ref{PMP1})-(\ref{PMP2}).
Following \cite{CaFr14, PAC} we shall adopt
\begin{definition}\label{def:2} For the optimization problem 
(\ref{Jbx})-(\ref{dyn2}) a point $\bar x \in \R^n$ is a {\bf  conjugate
point} if there exists $\ov z\in\R^n$ such that
$\bar x = x(0,\ov z)$, the control
$u(\cdot,\ov z)$ is a weak local minimizer of (\ref{Jbx})-(\ref{dyn2}), and  moreover
\bel{Conj-def}
\det\big(x_z(0,\ov z)\big)~=~0.
\eeq
Here $x_z$ denotes the $n\times n$ Jacobian matrix of partial derivatives of the 
map $z\mapsto x(0,z)$.
\end{definition}

Our main goal is to understand the structure of the set of conjugate points, for a 
generic terminal cost  $ \psi\in \C^5(\R^n)$
in (\ref{Jbx}).  The present paper provides two results in this direction.
In Section~\ref{sec:2} we prove a necessary condition for the optimality
of a trajectory starting at a conjugate point.
We recall that, by classical results \cite{Cesari, Clarke}, a trajectory 
$t\mapsto x(t)$ is not optimal if it contains a conjugate point $x(\tau)$ for some 
$0<\tau<T$.  However, the case $\tau=0$ is more delicate. A necessary condition that 
covers this case is given in Theorem~\ref{t:21}.
Relying on this more precise result, in Section~\ref{sec:3}
we study a classical problem in the Calculus of Variations:
$$\hbox{Minimize:}\quad \int_0^T L\bigl(\dot x(t)\bigr)\, dt + \psi\bigl(x(T)\bigr)
\qquad\hbox{subject to}~~x(0)=\bar x.$$
Assuming that the Lagrangian function $L=L(u)$ is smooth and uniformly convex,
 we study the structure of the set of conjugate points, for a generic terminal cost $\psi\in \C^4(\R^n)$.  In particular, we show that its $(n-2)$-dimensional
Hausdorff measure is locally finite.   In the 1-dimensional case, the set of conjugate points is empty.

\section{Necessary conditions for conjugate points}
\label{sec:2}
\setcounter{equation}{0}
In this section we derive  a necessary condition for conjugate points. 
For a given  $\ov z\in \R^n$, we consider the map $z\mapsto g(z,\ov z)$, defined by  
\bel{g}
g(z,\ov z)~\doteq~\int_0^{T}L\bigl(\Tilde{x}(t,z),u(t,z)\bigr)dt+\psi\bigl(\Tilde x (T,z)\bigr),
\eeq
where $u(t,z)\doteq u\bigl(x(t,z),p(t,z)\bigr)$ is the control corresponding to the solution of 
the backward Cauchy problem (\ref{PMP1})--(\ref{td3}), while $\Tilde{x}(\cdot,z)$ is the solution of 
\bel{fw}
\dot{x}(t)~=~f(x(t),u(t,z)),\qquad x(0)~=~x(0,\ov z).
\eeq
In other words,
$g(z,\ov z) $ is the cost of the trajectory $\Tilde x(\cdot,z)$ which
\begi
\item[(i)] starts  at the initial point $x(0, \ov z)$ of  the solution to 
the Pontryagin equations (\ref{PMP1})  ending at $\ov z$,
\item[(ii)] but uses the  control $u(\cdot, z)$, corresponding to the solution of (\ref{PMP1}) ending at $z$.
\endi

\begin{lemma}\label{Con-1} Let $\bar{z}\in \R^n$ and ${\bf v}\in \R^n$ be a 
unit vector such that $x_z(0,\ov z)\bfv=0$. Then the map $g_{{\bf v}}:\R\to\R$, defined by 
\[
g_{{\bf v}}(\theta)~\doteq~g(\ov z+\theta{\bf v},\ov z),
\] 
has first and second derivatives which vanish at $\theta=0$:
\bel{g'0}
g'_{{\bf v}}(0)~=~0,\qquad\qquad g''_{{\bf v}}(0)~=~0.
\eeq
\end{lemma}
{\bf Proof.} {\bf 1.} For a given solution to (\ref{PMP1})--(\ref{td3}), we denote by
$$x_z, p_z:[0,T]\mapsto \R^{n\times n},\qquad u_z:[0,T]\mapsto\R^{m\times n},$$
the matrix representations of the differentials w.r.t.~the terminal point $z$.
 Differentiating (\ref{PMP1}) one obtains
\bel{xp-y}\left\{
\bega{rl} 
\ds{d\over dt} x_z(t,z)&=~f_x(x,u)x_z+f_{u}(x,u) u_z\,,\\[4mm]
 \ds{d\over dt} \Tilde{x}_z(t,z)&=~f_x(\Tilde{x},u)\Tilde{x}_z+f_{u}(\Tilde{x},u)u_z\,,\\[4mm]
\ds{d\over dt} p_z(t,z)&=~-p_zf_x-p\left(f_{xx} x_z+f_{xu}u_z\right)-L_{xx}x_z-L_{xu}u_z\,.
\enda\right.
\eeq
Moreover, set  $\Gamma(t,z)\doteq L(x(t,z),u(t,z))$ for all $(t,z)\in [0,T]\times \R^n$. By (\ref{op1}) we have
\bel{k1}\bega{l}\ds
{d\over dt}\bigl[p(t,z)\cdot x_z(t,z)\bigr]~=~\bigl[ - p f_x - L_x\bigr] x_z + p \bigl[ f_x x_z+ f_u u_z\bigr]\\[4mm]
\qquad\qquad\ds
=~
-L_x x_z-L_{u}u_z~=~-\Gamma_z(t,z).\enda
\eeq
Observing that 
\bel{x=}
x(\cdot,\ov z)~=~\Tilde{x}(\cdot,\ov z),\qquad x_z(0,\ov z)\,\bfv~=~0,
\eeq
we have 
\bel{dz=eq}
\Tilde{x}_z(t,\ov z){\bf v}~=~x_z(t,\ov z){\bf v},\qquad \forall t\in [0,T].
\eeq
Recalling (\ref{g}), we now compute 
\bel{dg}
\begin{split}
g'_{\bf v}(\theta)&=~\int_{0}^{T}{d\over d\theta}L\bigl(\Tilde{x}(t,\ov z+\theta\bfv),u(t,\ov z+\theta{\bf v})\bigr)\,dt+\nabla\psi\bigl(\Tilde{x}(T,\bar{z}+\theta{\bf v})\bigr)\cdot {d\over d\theta} \Tilde{x}(T,\bar{z}+\theta{\bf v})\\
&=~\int_{0}^{T}\Gamma_z(t,\bar{z}+\theta{\bf v}){\bf v}\,dt+\nabla\psi\bigl(\Tilde{x}(T,\bar{z}+\theta{\bf v})\bigr)\cdot \Tilde x_z(T,
\bar{z}+\theta{\bf v})\bfv.
\end{split}
\eeq
Therefore, (\ref{k1})-(\ref{dz=eq}) yield
\[
\begin{split}
g'_{\bf v}(0)&~=~\int_{0}^{T}\Gamma_z(t,\bar{z}){\bf v}dt+\nabla\psi(\bar{z})\cdot
x_z(T,\ov z)\bfv\\
&~=~-\int_{0}^{T}{d\over dt}\big[p(t,\ov z)x_z(t,\ov z)\big]\bfv\,dt+p(T,\ov z)\,x_z(T,\ov z)\bfv~=~p(0,\ov z)\,x_z(0,\ov z)\bfv~=~0.
\end{split}
\]

{\bf 2.} To prove the second identity in (\ref{g'0}) one needs to differentiate 
(\ref{dg}) once more.   In the following, second order differentials such as  $\psi_{zz}=D^2_z \psi$  
and $\Tilde{x}_{zz}$ are regarded as 
symmetric bilinear maps, 
sending a couple of vectors $\bfv_1\otimes\bfv_2\in \R^n\times\R^n$ into $\R$ and into $\R^n$, respectively. 
We compute
\bel{g''}
\bega{rl}
g''_{\bf v}(\theta)&\ds=~\int_{0}^{T}{d^2\over d\theta^2}L(\Tilde{x}(t,\ov z+\theta\bfv),u(t,\ov z+\theta\bfv)\bigr) dt\\[3mm]
&\ds+ \psi_{zz}\bigl(\Tilde{x}(T,\ov z+\theta\bfv)\bigr)(\Tilde{\bfw}^{\theta}(T) \otimes\Tilde{\bfw}^{\theta}(T)) 
+\nabla\psi\bigl(\Tilde{x}(T,\ov z+\theta\bfv)
\bigr)\cdot \Big(\Tilde x_{zz}(T,\ov z+\theta\bfv)(\bfv\otimes \bfv)\Big),
\enda
\eeq
where
$$\Tilde{\bfw}^{\theta}(t)~\doteq~\Tilde{x}_z(t,\ov z+\theta\bfv)\bfv,\qquad t\in [0,T].$$

Differentiating the first equation in (\ref{xp-y}) once again w.r.t.~$z$, one obtains
\bel{d2x} {d\over dt} x_{zz}~=~f_{xx}(x_z\otimes x_z) +2 f_{xu}(x_z\otimes u_z) + f_{uu} (u_z\otimes u_z) 
+ f_x \,x_{zz} + f_u \,u_{zz}\,.\eeq
Setting
$$
{\bf w}^{\theta}(t)~\doteq~x_z(t,\ov z+\theta\bfv)\bfv,\qquad \bfb^\theta (t)~\doteq~u_z(t,\ov z+\theta\bfv)\bfv,\qquad t\in [0,T],
$$
{}from (\ref{d2x}) it follows
\bel{der2}
\bega{l} \ds
{d\over dt}x_{zz}(t,\bar{z})({\bf v}\otimes{\bf v})~=~f_{xx}\bigl(x(t,\bar{z}),u(t,\bar{z})\bigr)\bigl({\bf w}^{0}(t)\otimes {\bf w}^{0}(t)
\bigr)\\[4mm]
\qquad+2f_{xu}\bigl(x(t,\bar{z}),u(t,\bar{z})\bigr) \bigl(\bfw^{0}(t)\otimes \bfb^{0}(t)\bigr) +f_{uu}\bigl(x(t,\bar{z}),u(t,\bar{z})\bigr)\bigl(\bfb^{0}(t)\otimes \bfb^{0}(t)\bigr)\\[4mm]
\qquad+ 
f_x\bigl(x(t,\bar z),u(t,\bar{z})\bigr)x_{zz}(t,\bar{z})({\bf v}\otimes {\bf v})
+f_u\bigl(x(t,\bar{z}),u(t,\bar{z})\bigr)u_{zz}(t,\bar{z})({\bf v}\otimes {\bf v}),
\enda
\eeq
\bel{Tder2}
\bega{l} \ds
{d\over dt}\Tilde x_{zz}(t,\bar{z})({\bf v}\otimes{\bf v})~=~f_{xx}\bigl(\Tilde x(t,\bar{z}),u(t,\bar{z})\bigr)\bigl(\Tilde \bfw^{0}(t)\otimes \Tilde \bfw^{0}(t)
\bigr)\\[4mm]
\qquad+2f_{xu}\bigl(\Tilde x(t,\bar{z}),u(t,\bar{z})\bigr) \bigl(\Tilde \bfw^{0}(t)\otimes \bfb^{0}(t)\bigr)+f_{uu}\bigl(\Tilde x(t,\bar{z}),u(t,\bar{z})\bigr)\bigl(\bfb^{0}(t)\otimes \bfb^{0}(t)\bigr)\\[4mm]
\qquad+ 
f_x\bigl(\Tilde x(t,\bar{z}),u(t,\bar{z})\bigr)\Tilde x_{zz}(t,\bar{z})({\bf v}\otimes {\bf v})+f_u\bigl(\Tilde x(t,\bar{z}),u(t,\bar{z})\bigr)u_{zz}(t,\bar{z})({\bf v}\otimes {\bf v}).
\enda
\eeq
By (\ref{dz=eq}) one has
$$\bfw^0(t)~=~\Tilde \bfw^0(t)\qquad\qquad\forall t\in [0,T].$$
Comparing the two equations (\ref{der2})-(\ref{Tder2}), we see that  by (\ref{x=}) 
the only difference between the right hand sides
is the term involving $x_{zz}$.  Therefore we can write
\bel{txzz}
\Tilde x_{zz}(t,\ov z)({\bf v}\otimes{\bf v})~=~x_{zz}(t,\ov z)({\bf v}\otimes{\bf v})+w(t,\ov z),
\eeq
where $w(\cdot,\ov z):[0,T]\mapsto\R^n$ is the solution to the linear ODE
\bel{w}
\dot{w}(t)~=~f_x\bigl(x(t,\ov z),u(t,\ov z)\bigr)\cdot w(t),\qquad\qquad w(0)~=~-x_{zz}(0,\ov z)({\bf v}\otimes{\bf v}).
\eeq
Using (\ref{txzz}), we now compute
\bel{d2L}\bega{l}\ds
\left[{d^2\over d\theta^2}L\bigl(\Tilde{x}(t,\ov z+\theta\bfv),u(t,\ov z+\theta\bfv)\bigr)\right]_{\theta=0}\\[4mm]
\qquad\ds =~\Gamma_{zz}(t,\bar z)({\bf v}\otimes {\bf v})+
L_x\bigl(x(t,\bar{z}),u(t,\bar{z})\bigr)\,w(t,\bar{z}).
\enda\eeq
By (\ref{w}) and the second equation in (\ref{PMP1}) it follows
$${d\over dt}\bigl[p(t,\bar{z})w(t,\bar{z})\bigr] ~=~\dot p\, w + p f_x w 
~=~\bigl( - p f_x - L_x + pf_x \bigr) w~=~- L_x w\,.$$
Hence
\bel{Lxw}
L_x(x(t,\bar{z}),u(t,\bar{z}))\,w(t,\bar{z})~=~- {d\over dt}\bigl[p(t,\bar{z})w(t,\bar{z})\bigr].
\eeq
{}From  (\ref{g''}), using (\ref{k1}) and (\ref{Lxw}), and recalling that $\bfw^0(T)=\bfv$ while 
$\psi_z(\bar z) = p(T, \bar z)$, we obtain
\[
\bega{rl}
g''_{\bf v}(0)&\ds=~\int_{0}^{T}\Gamma_{zz}(t,\bar z)({\bf v}\otimes {\bf v})\,dt+\int_{0}^{T}L_x(x(t,\bar{z}),u(t,\bar{z}))w(t,\bar{z})dt\\[4mm]
&\qquad\quad \ds +\psi_{zz}(\ov z)(\bfv \otimes \bfv)+\psi_z(\ov z)\cdot \Big[x_{zz}(T,\ov z)(\bfv\otimes \bfv)+w(T,\bar{z})\Big]\\[4mm]
&\ds=~-\int_{0}^{T}{d\over dt}\left({d\over dz}\bigl[p(t,\bar{z})\, x_z(t,\bar{z})\bigr]\right)({\bf v}\otimes {\bf v})\,dt
-\int_{0}^{T}{d\over dt}\left[p(t,\bar{z})w(t,\bar{z})\right]dt\\[4mm]
&\qquad\quad\ds  +{d\over dz}[p(T,\bar{z})x_z(T,\bar{z})]({\bf v}\otimes{\bf v})+p(T,\bar{z})w(T,\bar{z})\\[4mm]
&\ds=~{d\over dz}\bigl[p(0,\bar{z})\cdot x_z(0,\bar{z})\bigr]({\bf v}\otimes {\bf v})+p(0,\bar{z})w(0,\bar{z})\\[4mm]&=~p(0,\bar{z})\bigl[x_{zz}(0,\bar{z})({\bf v}\otimes{\bf v})+w(0,\bar{z})\bigr]+
\bigl(p_z(0,\bar{z})\, \bfv\bigr)\cdot\bigl(x_z(0,\bar{z})\,\bfv\bigr)\\[4mm]
&=~p(0,\bar z) \, \Tilde x_{zz}(0,\bar z)({\bf v}\otimes{\bf v})~=~0.
\enda
\]
The proof is complete.
\endproof

In view of (\ref{g'0}), if $\bar x = x(0, \bar z)$ is a conjugate point the optimality assumption implies
the vanishing of the third derivative:
\bel{g30}
g_\bfv'''(0)~=~0.\eeq
This yields the following necessary condition:

\begin{theorem}\label{t:21} 
Given  a conjugate point $\bar{x}=x(0,\ov z)\in \R^n$, with $\bar{z}\in\R^n$ associated to   a weak local minimizer $u(\cdot,\ov z)$  of the optimization problem (\ref{Jbx})-(\ref{dyn2}),  let ${\bf v}\in \R^n$ be a 
unit vector such that $x_z(0,\ov z)\bfv=0$. Then  one has
\bel{ke}
\bigl(p_z(0,\ov z)\bfv\bigr)\cdot x_{zz}(0,\ov z)({\bf v}\otimes{\bf v})~=~0.
\eeq
\end{theorem}
{\bf Proof.} Differentiating (\ref{g''}), we compute 
\bel{g'''}
\bega{rl}
g'''_{\bf v}(\theta)&\ds=~\int_{0}^{T}{d^3\over d\theta^3}L(\Tilde{x}(t,\ov z+\theta\bfv),u(t,\ov z+\theta\bfv)\bigr)\,dt\\[4mm]
&\qquad \ds+\psi_{zzz}\bigl(\Tilde{x}(T,\ov z+\theta\bfv)\bigr)\bigl(\Tilde{\bfw}^{\theta}(T) \otimes\Tilde{\bfw}^{\theta}(T)\otimes \Tilde{\bfw}^{\theta}(T) \bigl)  \\[3mm]
&\qquad +3\psi_{zz}\bigl(\Tilde{x}(T,\ov z+\theta\bfv)\bigr)\bigl(\Tilde{\bfw}^{\theta}(T) \otimes \Tilde x_{zz}(T,\ov z+\theta\bfv)(\bfv\otimes {\bf v})\bigr)\\[3mm]
&\qquad\ds + \psi_z\bigl(\Tilde{x}(T,\ov z+\theta\bfv)
\bigr)\cdot \Big(\Tilde x_{zzz}(T,\ov z+\theta\bfv)(\bfv\otimes \bfv\otimes {\bf v})\Big).
\enda
\eeq
The third order differentials  $ \psi_{zzz}$  and $\Tilde{x}_{zzz}$ are here regarded  as 
tri-linear maps, 
sending a triple of vectors $\bfv_1\otimes\bfv_2\otimes {\bf v}_3\in \R^n\times\R^n\times \R^n$ into $\R$ and  into $\R^n$, respectively.  Differentiating the identities (\ref{der2})-(\ref{Tder2})  once more w.r.t.~$z$, we obtain
\bel{der3}
\bega{l} \ds
{d\over dt}x_{zzz}(t,\bar{z})({\bf v}\otimes{\bf v}\otimes{\bf v})~=~f_{xxx}\bigl(x(t,\bar{z}),u(t,\bar{z})\bigr)\bigl({\bf W}_0(t)\bigr)+3f_{xxu}\bigl(x(t,\bar{z}),u(t,\bar{z})\bigr)\bigl({\bf W}_1(t)\bigr)\\[4mm]
\qquad +3f_{xuu}\bigl(x(t,\bar{z}),u(t,\bar{z})\bigr) \bigl({\bf W}_2(t)\bigr)+f_{uuu}\bigl(x(t,\bar{z}),u(t,\bar{z})\bigr)\bigl({\bf W}_3(t)\big)\\[4mm]
\qquad+3f_{xx}\bigl(x(t,\bar{z}),u(t,\bar{z})\bigr)\bigl(\bfw^{0}(t)\otimes {\bf x}(t)\bigr)+3f_{uu}\bigl(x(t,\bar{z}),u(t,\bar{z})\bigr)\bigl({\bf u}(t)\otimes {\bf b}^0(t)\bigr)\\[4mm]
\qquad +3f_{xu}\bigl(x(t,\bar{z}),u(t,\bar{z})\bigr)\bigl( {\bf x}(t)\otimes {\bf b}^0(t)+{\bf u}(t)\otimes \bfw^{0}(t)\bigr)\\[4mm]
\qquad+f_u\bigl(x(t,\bar{z}),u(t,\bar{z})\bigr){\bf U}(t)+f_x\bigl(x(t,\bar{z}),u(t,\bar{z})\bigr){\bf X}(t)
\enda
\eeq
and 
\bel{Tder3}
\bega{l} \ds
{d\over dt}\Tilde{x}_{zzz}(t,\bar{z})({\bf v}\otimes{\bf v}\otimes{\bf v})~=~f_{xxx}\bigl(x(t,\bar{z}),u(t,\bar{z})\bigr)\bigl({\bf W}_0(t)\bigr)+3f_{xxu}\bigl(x(t,\bar{z}),u(t,\bar{z})\bigr)\bigl({\bf W}_1(t)\bigr)\\[4mm]
\qquad +3f_{xuu}\bigl(x(t,\bar{z}),u(t,\bar{z})\bigr) \bigl({\bf W}_2(t)\bigr)+f_{uuu}\bigl(x(t,\bar{z}),u(t,\bar{z})\bigr)\bigl({\bf W}_3(t)\big)\\[4mm]
\qquad+3f_{xx}\bigl(x(t,\bar{z}),u(t,\bar{z})\bigr)\bigl(\bfw^{0}(t)\otimes \Tilde{{\bf x}}(t)\bigr)+3f_{uu}\bigl(x(t,\bar{z}),u(t,\bar{z})\bigr)\bigl({\bf u}(t)\otimes {\bf b}^0(t)\bigr)\\[4mm]
\qquad+3f_{xu}\bigl(x(t,\bar{z}),u(t,\bar{z})\bigr)\bigl( \Tilde{{\bf x}}(t)\otimes {\bf b}^0(t)+{\bf u}(t)\otimes \bfw^{0}(t)\bigr)\\[4mm]
\qquad+f_u\bigl(x(t,\bar{z}),u(t,\bar{z})\bigr){\bf U}(t)+f_x\bigl(x(t,\bar{z}),u(t,\bar{z})\bigr)\Tilde{{\bf X}}(t),
\enda
\eeq
with 
\[
\begin{cases}
{\bf W}_0(t)~\doteq~{\bf w}^{0}(t)\otimes {\bf w}^{0}(t)\otimes{\bf w}^{0}(t),\qquad {\bf W}_1(t)~\doteq~{\bf w}^{0}(t)\otimes {\bf w}^{0}(t)\otimes{\bf b}^{0}(t),\\[3mm]
{\bf W}_2(t)~\doteq~{\bf w}^{0}(t)\otimes {\bf b}^{0}(t)\otimes{\bf b}^{0}(t),\qquad {\bf W}_3(t)~\doteq~{\bf b}^{0}(t)\otimes {\bf b}^{0}(t)\otimes{\bf b}^{0}(t),\\[3mm]
{\bf X}(t)~\doteq~x_{zzz}(t,\bar{z})({\bf v}\otimes {\bf v}\otimes {\bf v}),\quad \Tilde{{\bf X}}(t)~\doteq~\Tilde{x}_{zzz}(t,\bar{z})({\bf v}\otimes {\bf v}\otimes {\bf v}),\quad{\bf x}(t)~\doteq~x_{zz}(t,\bar{z})({\bf v}\otimes {\bf v}) \\[3mm]
\Tilde{{\bf x}}(t)~\doteq~\Tilde{x}_{zz}(t,\bar{z})({\bf v}\otimes {\bf v}),\quad {\bf u}(t)~\doteq~u_{zz}(t,\bar{z})({\bf v}\otimes {\bf v}),\quad {\bf U}(t)~\doteq~u_{zzz}(t,\bar{z})({\bf v}\otimes {\bf v}\otimes {\bf v})
\end{cases}
\]
Comparing
the results, we eventually obtain
\bel{Txzzz}
\Tilde x_{zzz}(t,\ov z)({\bf v}\otimes{\bf v}\otimes{\bf v})~=~x_{zzz}
(t,\ov z)({\bf v}\otimes{\bf v}\otimes{\bf v})+W(t,\ov z),
\eeq
where $w(\cdot)$ is the function constructed at (\ref{w}), while $W(\cdot,\ov z):[0,T]\to \R^n$ is the solution to the linear ODE
\bel{3ode}
\bega{rl}
\dot{W}(t)&=~f_x\bigl(x(t,\bar{z}),u(t,\bar{z})\bigr)\cdot W(t)+3f_{xx}\bigl(x(t,\bar{z}),u(t,\bar{z})\bigr)\bigl ({\bf w}^0(t)\otimes w(t,\bar{z})
\bigr)\\[3mm]
&\qquad\qquad\qquad\qquad +3 f_{xu}(x(t,\bar{z}),u(t,\bar{z}))\bigl(w(t,\bar{z})\otimes {\bf b}^0(t)\bigr),
\enda\eeq
with initial data
\bel{W0}
 W(0,\bar{z})~=~-x_{zzz}(0,\ov z)({\bf v}\otimes{\bf v}\otimes{\bf v}).
\eeq
%
%

In this case, we have 
\[
\begin{split}
\Big[{d^3\over d\theta^3}&L\bigl(\Tilde{x}(t,\ov z+\theta\bfv),u(t,\ov z+\theta\bfv)\bigr)\Big]_{\theta=0}\\
&~=~\Gamma_{zzz}(t,\bar z)({\bf v}\otimes{\bf v}\otimes {\bf v})
+L_x\,W(t,\ov z)\\
&\qquad\qquad+3L_{xu}\,\bigl(w(t,\ov z)\otimes {\bf b}^0(t)\bigr)
+3L_{xx}\,\bigl(w(t,\ov z)\otimes {\bf w}^{0}(t)\bigr),
\end{split}
\]
\[
\begin{split}
{d\over dt}\Big[p(t,\bar{z})\,W(t,\bar{z})\Big]&~=\,- L_x\,W(t,\bar{x})
+3 p(t,\bar{z})f_{xx}\,\bigl({\bf w}^0(t)\otimes w(t,\bar{z})\bigr)\\
&\qquad\qquad\qquad +3 p(t,\bar{z})f_{xu}\,\bigl( w(t,\bar{z})\otimes {\bf b}^0(t)
\bigr),
\end{split}
\]
\[
\bega{rl}\ds
{d\over dt}\Big[\bigl(p_z(t,\ov z)\bfv \bigr) \cdot w(t,\bar{z})\Big]&\ds=\,-p(t,\bar{z})f_{xx}\,
\bigl({\bf w}^0(t)\otimes w(t,\bar{z})\bigr)-p(t,\bar{z})f_{xu}\, \bigl({\bf b}^0(t)\otimes w(t,\bar{z})\bigr)\bigr)\\[3mm]
&\ds\qquad 
-L_{xx}\,\bigl({\bf w}^0(t)\otimes w(t,\bar{z})\bigr)
 -L_{xu}\,\bigl( w(t,\bar{z})\otimes  {\bf b}^0(t)\bigr).
\enda
\]
In the above formulas, it is understood that the functions $f, L$ and all their partial derivatives are 
computed 
at the point $\bigl(x(t,\bar{z}),u(t,\bar{z})\bigr)$.

Using the above identities together with (\ref{k1}) and  (\ref{w}), from (\ref{g'''}) we obtain
{\small
\[
\bega{l}\ds
g'''_{{\bf v}}(0)~=~\left[-\int_{0}^{T}{d\over dt}\left({d^2\over dz^2}\bigl[p(t,\bar{z})\cdot x_z(t,\bar{z})\bigr]({\bf v}\otimes{\bf v}\otimes {\bf v})\right)\,dt+{d^2\over dz^2}\bigl[p(T,\bar{z})\cdot x_z(T,\bar z)\bigr]({\bf v}\otimes{\bf v}\otimes {\bf v})\right]\\[4mm]
\ds  -\left[\int_{0}^{T}{d\over dt}\left[pW\right](t,\bar{z})dt-\left[pW\right](T,\bar{z})\right]+3\left[\bigl(p_z(T,\ov z){\bf v}\bigr)\cdot w(T,\bar{z})-\int_{0}^{T}{d\over dt}\big[\bigl(p_z(t,\ov z){\bf v}\bigr)\cdot w(t,\bar{z})\big]\,dt\right]\\[4mm]
\quad =~\ds {d^2\over dz^2}\bigl[p(0,\bar{z})\cdot x_z(0,\bar z)\bigr]({\bf v}\otimes{\bf v}\otimes {\bf v})+[pW](0,\bar{z})
+3 \bigl(p_z(0,\ov z)\,\bfv\bigr) \cdot w(0,\bar{z})\\[3mm]
\quad =~2 \bigl(p_z(0,\ov z)\,\bfv\bigr) \cdot
x_{zz}(0,\bar{z})({\bf v}\otimes {\bf v})+3 \bigl(p_z(0,\ov z)\,\bfv\bigr) \cdot w(0,\bar{z})\\[3mm]
\quad =~-\bigl(p_z(0,\ov z)\,\bfv\bigr) \cdot
x_{zz}(0,\bar{z})({\bf v}\otimes {\bf v}).
\enda
\]
}
Since $g_\bfv$ attains a local minimum at $\theta=0$ and $g'_{{\bf v}}(0)=g''_{{\bf v}}(0)=0$,
this yields  (\ref{ke}).
\endproof
\section{Conjugate points for a generic problem in the  Calculus of Variations}
\label{sec:3}
\setcounter{equation}{0}
In this section, the necessary condition stated in Theorem~\ref{t:21} will be used
to study a generic property of the set of conjugate points for a classical problem in the 
Calculus of Variations. Namely, we seek to minimize (\ref{Jbx}) in the special case
where
\bel{Ldef}
\dot{x}~=~u,\qquad L(x,u)~=~L(u).
\eeq
In this case (see for example \cite{BPi}), the value function $V$ is the unique viscosity solution to the Hamilton-Jacobi equation 
\bel{HJ}
\begin{cases}
-V_{t}(t,x)- H\bigl(\nabla V(t,x)\bigr)~=~0,\qquad (t,x)\in [0,T]\times \R^n,\\[4mm]
V(T,x)~=~\psi(x),\qquad\qquad\qquad  x\in \R^n,
\end{cases}
\eeq
with 
\bel{H-c}
H(p)~=~\min_{\omega\in\R^n}\bigl\{L(\omega)+p\cdot\omega\bigr\}.
\eeq
By (\ref{PMP1}) and (\ref{td3}) it follows
\bel{zp}
x(0,z)~=~z- T\cdot DH(\nabla \psi(z)),\qquad p(0,z)~=~\nabla\psi(z),\qquad z\in \R^n.
\eeq
By Theorem~\ref{t:21}, the conjugate points are thus contained in the set 
$\big\{x(0,z)\,;~ (z,\bar{v})\in\Omega_{\psi}\big\} $, with
\bel{J}  \Omega_{\psi}~\doteq~\Big\{(z,\bfv)\in\R^n\times S^{n-1}\,;~x_z(0, z)\bfv =0, \quad \bigl(D^2\psi(z)\bfv\bigr)\cdot x_{zz}(0,z)({\bf v}\otimes {\bf v})=0\Big\},
\eeq
where $S^{n-1}$ denotes the set of unit vectors in $\R^n$.

\begin{theorem}\label{Cj}  Let the function $L=L(u)$ be smooth and uniformly convex. 
Then there exists a $\mathcal{G}_{\delta}$ subset $\mathcal{M}\subseteq\mathcal{C}^4(\R^n)$ such that for every $\psi\in \mathcal{M}$, the set $\Omega_{\psi}$  at (\ref{J})
  is  
   an embedded manifold  of dimension $n-2$. 
  

  \end{theorem}
  
{\bf Proof.} {\bf 1.} Given a terminal cost $\psi\in\mathcal{C}^4(\R^n)$,  defining  $(x,p)$ as in  (\ref{zp}), 
we have
\bel{37}
x_z(0,z)~=~\mathbb{I}_n-T\cdot D^2H(\nabla\psi(z))D^2\psi(z),\qquad p_z(0,z)~=~D^2\psi(z).
\eeq
Thus, if $x_z(0,z){\bf v}=0$ then ${\bf v}=T\cdot D^2H(\nabla\psi(z))D^2\psi(z){\bf v}$ and 
\bel{38}
p_z(0,z){\bf v}~=~D^2\psi(z){\bf v}~=~{1\over T}\cdot \big[D^2H(\nabla\psi(z))\big]^{-1}({\bf v}).
\eeq
This implies 
\bel{Omega-psi}
 \Omega_{\psi}~=~\left\{(z,{\bf v})\in\R^n\times S^{n-1}\,;~~x_z(0, z)\bfv =0, ~~ \big[D^2H(\nabla\psi(z))\big]^{-1}{\bf v}\cdot x_{zz}(0,{z})({\bf v}\otimes {\bf v})=0
\right\}.
\eeq
Define the $\C^1$ map  $\Phi^{\psi}:\R^n\times S^{n-1}\mapsto\R^n\times\R$ 
by setting
\bel{Phi-psi}
\Phi^{\psi}(z,{\bf v})~\doteq~\left(x_z(0,z){\bf v} \,, ~ \big[D^2H(\nabla\psi(z))\big]^{-1}{\bf v}\cdot  x_{zz}(0,z)({\bf v}\otimes {\bf v})\right).
\eeq
For $k\geq 1$, let $\overline{B}_k\subset\R^n$ be the closed ball centered at the origin with radius $k$, and consider the open subset of $\mathcal{C}^4(\R^n)$
\bel{M-k}
\mathcal{M}_{k}~\doteq~\left\{\psi\in \mathcal{C}^4(\R^n) \,:~~ \Phi^{\psi}\vert_{\overline{B}_k\times S^{n-1}}\text{ is transversal to } \{\bf0\}\right\}.
\eeq
Here $\{\bf0\}$ denotes 
the zero-dimensional manifold containing the single point $(0,0)\in \R^n\times\R$.
%

If $\mathcal{M}_k$ is dense in $\mathcal{C}^4(\R^n)$ for all $k\in\mathbb{Z}^+$, then the set  $\ds\mathcal{M}\doteq \bigcap_{k\geq 1} \mathcal{M}_k$ is a $\mathcal{G}_{\delta}$ subset of $\mathcal{C}^4(\R^n)$ such that for every $\psi\in \mathcal{M}$, $\Phi^{\psi}$ is transverse to $\{\bf 0\}$. By the implicit function theorem, the set $\Omega_{\psi}$ is an embedded manifold  of dimension $n-2$.

%
\v

{\bf 2.} Next, we show that $\M_k$ is dense in $\C^4(\R^n)$. 
For this purpose, fix any $\Tilde{\psi}\in \C^4(\R^n)$.  For every $\ve>0$, we first approximate $\Tilde{\psi}$ by a smooth function $\psi$ with $\|\Tilde{\psi}-\psi\|_{\C^4}<\ve$. Then we need to construct a perturbed 
function $\psi^\theta$ arbitrarily close to $\psi$ in the $\C^4$ norm, which 
lies in $\M_k$.  Toward this goal, for any point $(\bar z, \bfv)\in \ov B_k\times S^{n-1}$,
we consider the family of perturbed functions of the form
\bel{psite}\psi^\theta(z)~\doteq~\psi(z) +\eta(z-\bar z)\cdot
\left[ \sum_{i,j=1}^n \theta_{ij} (z_i-\bar z_i) (z_j-\bar z_j)
+ \sum_{k=1}^n \theta_k (z_k-\bar z_k)^3\right].\eeq
Here $\eta:\R^n\mapsto [0,1]$ is a smooth cutoff function, such that 
\bel{cutoff}\eta(y)~=~\left\{ \bega{rl} 1\quad &\hbox{if}~~|y|\leq 1,\\[2mm]
0\quad &\hbox{if}~~|y|\geq 2.\enda\right.\eeq
Moreover, $\theta = (\theta_{ij}, \theta_k)\in \R^{n^2+n}$.
We claim that the map 
$$(z,\bfv, \theta)~\mapsto~\Phi^{\psi^\theta}~=~ \Big(x^{\theta}_z(0,z){\bf v} , \   \big[D^2H(\nabla\psi^\theta(z))\big]^{-1}{\bf v}\cdot x^{\theta}_{zz}(0,{z})({\bf v}\otimes {\bf v})\Big)$$
is transversal to $\{{\bf 0 }\}\subset\R^n\times\R$ at the point $(\bar z, \bfv, 0)$.
This will certainly be true if the Jacobian matrix $D_\theta \Phi^{\psi^\theta}$
of partial derivatives w.r.t.~$\theta_{ij}, \theta_k$ has maximum rank $n+1$.

Writing $z=(z_1,\ldots, z_n)\in \R^n$ and recalling (\ref{psite}),
for $|z-\bar z|<1$
we compute the partial derivatives 
\bel{1der}
{\partial\over \partial z_i}\psi^\theta(z)~=~{\partial\over \partial z_i}\psi(z)+ 
\sum_{j=1}^n (\theta_{ij}+\theta_{ji})  (z_j-\bar z_j) + 3  \theta_i (z_i-\bar z_i)^2,
\eeq
\bel{2der}
{\partial^{2}\over \partial z_i\partial z_j}\psi^\theta(z)~=~{\partial^{2}\over \partial z_i\partial z_j}\psi(z)+ \left\{\bega{cl} 2\theta_{ii}+ 6  \theta_i (z_i-\bar z_i)\quad &\hbox{if} \quad i=j,\\[2mm]
\theta_{ij}+\theta_{ji}\quad &\hbox{if} \quad i\not= j.\enda\right.
\eeq
Calling $\{{\bf e}_1,\cdots, {\bf e}_n\}$ the standard basis of $\R^n$, we have 
\[
D^2\psi^{\theta}(z){\bf v}~=~D^2\psi(z){\bf v}+\sum_{i=1}^{n}\sum_{j=1}^{n}(\theta_{ij}+\theta_{ji}){\bf v}_j\cdot {\bf e}_i+6\sum_{i=1}^{n}\theta_i (z_i-\bar{z}_i){\bf v}_i\cdot {\bf e}_i\,,
\]
\[
{\partial\over \partial\theta_{ij}}D^2\psi^{\theta}(z){\bf v}~=~{\bf v}_j\cdot {\bf e}_i+{\bf v}_i\cdot {\bf e}_j\,.
\]
 Thus, for every $i\in\{1,\cdots,n\}$,  the matrix $D_{\theta}[D^2\psi^{\theta}(\bar z){\bf v}]$ contains the $n\times n$  submatrix
$$S_i~\doteq~\left[\ds{\partial\over \partial\theta_{ij}}D^2\psi^{\theta}(\bar z){\bf v}\right]_{j=1}^{n}~=~\begin{bmatrix}{\bf r}_1\\ \vdots\\ {\bf r}_i\\   \vdots \\{\bf r}_n\end{bmatrix}
~=~
\begin{bmatrix} \bfv_i &0&\cdots &0&\cdots&\cdots&0\cr
0& \bfv_i&\cdots &0&\cdots&\cdots&0\cr
\vdots&\vdots & \ddots&\vdots &\cdots&\cdots&0\cr
{\bf v}_1& \cdots & {\bf v}_{i-1} &2{\bf v}_i &{\bf v}_{i+1}&\cdots&{\bf v}_n\cr
0&\cdots&0&0&\bfv_i&\cdots&0\cr
\vdots&&\vdots&\vdots&\vdots&\ddots & \vdots\cr
0&\cdots&0 &0& 0&\cdots &\bfv_i\end{bmatrix}.
$$
Notice that this implies
\bel{Si}
\det(S_i)~=~2{\bf v}_i^n.
\eeq
Recalling (\ref{37}) and (\ref{1der}),  we have 
\bel{con1}
D_{\theta}\big[x_z^{\theta}(0,\bar{z}){\bf v}\big]~=~-T\cdot D^2H(\nabla\psi^{\theta}(\bar{z}))D_{\theta}[D^2\psi^\theta(\bar{z}){\bf v}\big], 
\eeq
and  (\ref{Phi-psi})  implies that   the matrix $D_{\theta}\big[\Phi^{\psi^\theta}(\bar{z},\bfv,\theta)\big]$ contains the $(n+1)\times n$ submatrix 
\bel{con2}
\Tilde{S}_{i}~\doteq~-T\cdot\begin{bmatrix}D^2H(\nabla\psi^{\theta}(\bar{z}) )S_i\\[2mm]
0_{1\times n}
\end{bmatrix}.
\eeq
Furthermore, we compute
\[
x^{\theta}_{zz}(0,z)({\bf v}\otimes {\bf v})~=~-TD^2H(\nabla\psi^{\theta}(z))D^3\psi^{\theta}(z)({\bf v}\otimes {\bf v})-TD^3H(\nabla\psi^{\theta}(z))(D^2\psi^{\theta}(z){\bf v}\otimes D^2\psi^{\theta}(z){\bf v}),
\]
with
\[
\begin{cases}
\ds{\partial\over \partial\theta_i}D\psi^{\theta}(z)~=~3(z_{i}-\bar z_i)^2\cdot{\bf e}_i,\qquad \ds{\partial\over \partial\theta_i}D^2\psi^{\theta}(z){\bf v}~=~6(z_i-\bar z_i){\bf v}_i\cdot {\bf e}_i,\\[5mm]
\ds{\partial\over \partial\theta_i}D^3\psi^{\theta}(z)({\bf v}\otimes {\bf v})~=~6\bfv^2_i\cdot 
{\bf e}_i,\qquad\forall i\in \{1,\cdots, n\}.
\end{cases}
\]
In particular, we have
\[
\ds{\partial\over \partial\theta_k}x^{\theta}_{zz}(\bar z)({\bf v}\otimes {\bf v})~=~-6T\bfv_k^2D^2H(\nabla\psi^{\theta}(\bar z)){\bf e}_k\,,
\]
and this yields 
\bel{ke1}
\begin{split}
\ds{\partial\over \partial\theta_k} \big(\big[D^2H(\nabla\psi^{\theta}(\bar z))\big]^{-1}\bfv\cdot x^{\theta}_{zz}&(0,\bar{z})({\bf v}\otimes {\bf v})\big)~=~\big[D^2H(\nabla\psi^{\theta}(\bar z))\big]^{-1}\bfv\cdot \ds{\partial\over \partial\theta_k}x^{\theta}_{zz}(\bar z)({\bf v}\otimes {\bf v})\\
&~=~-6T {\bf v}^2_k\cdot \big[D^2H(\nabla\psi^{\theta}(\bar z))\big]^{-1}\bfv\cdot D^2H(\nabla\psi^{\theta}(\bar z)){\bf e}_k\\
&~=~-6T {\bf v}^2_k\cdot {\bf v}\cdot {\bf e}_k~=~-6T {\bf v}_k^3.
\end{split}
\eeq
By the previous analysis we conclude that, for every $i\in \{1,\cdots, n\}$, 
the Jacobian matrix
$D_\theta \Phi^{\psi^\theta}$
of partial derivatives w.r.t.~$\theta_{ij}, \theta_i$
contains $n+1$ columns which form the $(n+1)\times n$ submatrix
\bel{con3}
\Lambda_i~\doteq~[\Tilde{S}_i,b_i]\qquad\mathrm{with}\qquad b_{i}=(*\,,\,\cdots\,,\,*\,,\, -6T {\bf v}_i^3)^{T}.
\eeq
By  (\ref{Si}) and  (\ref{con2}), it follows
\[
\det(\Lambda_i)~=~- 6T {\bf v}_i^3\cdot\det\left(-T\cdot D^2H(\nabla\psi^{\theta}(\bar z) )S_i\right)~=~-12T{\bf v}_i^{n+3}\cdot \det\left(-T\cdot D^2H(\nabla\psi^{\theta}(\bar z) )\right).
\]
%
By the strict convexity of $H$ and since ${\bf v}\in S^{n-1}$, we have that $\mathrm{rank}(\Lambda_i)=n+1$ for some $i\in\{1,\cdots,n\}$ and this yields 
\bel{fullrank}
\mathrm{rank} \,D_\theta \Phi^{\psi^\theta}(\bar z,\bfv,0)~=~n+1.
\eeq

{\bf 3.} By continuity, there exists a  neighborhood ${\cal N}_{\bar z, \bfv}$ of $(\bar z, \bfv)$ such that
\[
 \mathrm{rank} \,D_\theta \Phi^{\psi^\theta}(\bar z',\bfv',0)~=~n+1\qquad\forall (\bar z', \bfv')\in {\cal N}_{\bar z, \bfv}.
\]
%
Covering the compact set $\ov B_k\times S^{n-1}$ 
with finitely many open neighborhoods ${\cal N}_\ell=
{\cal N}_{\bar z^\ell, \bfv^\ell}$, $\ell=1,\ldots,N$, we consider the family of combined perturbations
\bel{psiall}\psi^\theta(z)~=~\psi(z) +\sum_{\ell=1}^N\left\{ \eta(z-\bar z^\ell)\cdot
\left[ \sum_{i,j=1}^n \theta^\ell_{ij} (z_i-\bar z^\ell_i) (z_j-\bar z^\ell_j)
+ \sum_{k=1}^n \theta_k^\ell (z_k-\bar z^\ell_k)^3\right]\right\}.\eeq
By construction,  the matrix of partial derivatives
w.r.t.~all combined variables $\theta = (\theta^\ell_{ij}, \theta^\ell_k)$ satisfies 
(\ref{fullrank}) at every point $(\bar z,\bfv)\in \ov B_k\times S^{n-1}$.

Again by continuity, we still have
\bel{fullrank}
\mathrm{rank} \,D_\theta \Phi^{\psi^\theta}(\bar z,\bfv,\theta)~=~n+1
\eeq
for all $(\bar z,\bfv,\theta)\in \ov B_k\times S^{n-1}\times\R^{N(n^2+n)}$ with $|\theta|<\ve$
sufficiently small.   By the transversality theorem \cite{Bloom, GG},
this implies that, for a dense set of values $\theta$, the smooth map $\Phi^{\psi^\theta}$
at (\ref{Phi-psi})
is transversal to $\{\bf 0\}$, restricted to the domain $(\bar z, \bfv)\in \ov B_k\times S^{n-1}$. We conclude that the set $\mathcal{M}_k$ is dense and the proof is complete.
\endproof

\begin{corollary}\label{Cj2} As the same setting in Theorem \ref{Cj}, there exists a $\mathcal{G}_{\delta}$ subset $\mathcal{M}\subseteq\mathcal{C}^4(\R^n)$ with the following property.
For every $\psi\in \mathcal{M}$, the set  $\Gamma_\psi\subset\R^n$  of all conjugate points  has locally bounded $(n-2)$-dimensional Hausdorff measure.
\end{corollary}

{\bf Proof.} Call $\pi:\R^n\times S^{n-1}\to\R^n$   the projection on the first component, so that $\pi(z,{\bf v})=z$.  Then the set of all conjugate points satisfies
the inclusion
\[
\Gamma_{\psi}~\subseteq~ \Big\{x\bigl(0,\pi(z,{\bf v})\bigr)\,\doteq\, 
z-T\cdot DH\bigl(\nabla \psi(z)\bigr)\,;~~(z,{\bf v})\in\Omega_\psi\Big\}.
\]
By Theorem \ref{Cj}, there exists a $\G_\delta$ set $\M\subset\C^4(\R^n)$ such that, for $\psi\in \M$, the set 
 $\Omega_\psi$ is an 
 embedded manifold  of dimension $n-2$. 
 
 We now observe that the map $(z,{\bf v})\mapsto x(0,\pi(z,{\bf v}))$ is  Lipschitz continuous. Moreover, for every $z\in\R^n$, one has 
 \[ |z|~\leq~
\bigl|x(0,z)\bigr| +LT\,,\qquad\mathrm{with}\qquad L~\doteq~\max_{|p|\leq \|\nabla\psi\|_{\infty}}|DH(p)|,
\]
and this implies
$$\Gamma_\psi\cap \ov B_r
~\subseteq~\Big\{x(0,z)\,;~~(z,{\bf v})\in\Omega_\psi\,,~~z\in \ov B_{r+LT}\Big\}\qquad\forall r>0.
$$
Since $\Omega_\psi$ is an embedded manifold,  and the map $(z,\bfv)\mapsto
x(0,z)$ is Lipschitz  continuous, by the properties of Hausdorff measures \cite{EG} we conclude that  the set  $\Gamma_\psi$  has locally bounded $(n-2)$-dimensional Hausdorff measure.
\endproof

\begin{remark} {\rm 
Using the original version of Sard's theorem \cite{Sard}, the smoothness assumption on $L$ can be somewhat relaxed.    Indeed, one can check that 
both Theorem~\ref{Cj} and Corollary~\ref{Cj2} still hold for a uniformly convex Lagrangian function $L\in \mathcal{C}^{n+2}$.}
\end{remark}

{\bf Acknowledgment.} The research of  Khai T.~Nguyen  was partially supported by National Science Foundation with grant DMS-2154201.
%
%

\end{document}